\documentclass{amsart}
\usepackage{amssymb}
\usepackage{amsthm}
\usepackage{hyperref}
\usepackage{amsmath}
\usepackage{amssymb}
\usepackage{enumerate}
\usepackage{graphicx}
\usepackage{amsthm}
\usepackage{amsopn}
\usepackage{amsfonts}
\usepackage{upgreek}
\usepackage{amsfonts}
\usepackage{amssymb}
\usepackage{empheq}
\usepackage{caption}
\numberwithin{equation}{section}
\usepackage{comment}
\usepackage{lineno}

\newtheorem{theorem}{Theorem}

\newtheorem{example}{Example}%
\newtheorem{remark}{Remark}%
\newtheorem{lemma}{Lemma}%
\newtheorem{definition}{Definition}%

\begin{document}
\title[A simple oscillation criteria for second-order DEPCAG]{A simple oscillation criteria for second-order differential equations with piecewise constant argument of generalized type}
\author{Ricardo Torres Naranjo}
\address{Instituto de Ciencias Físicas y Matemáticas, Facultad de Ciencias, Universidad Austral de Chile\\
Campus Isla Teja s/n, Valdivia, Chile.}
\address{\noindent Instituto de Ciencias, Universidad Nacional de General Sarmiento\\
Los Polvorines, Buenos Aires, Argentina}
\curraddr{}
\email{ricardo.torres@uach.cl}
\thanks{}
\subjclass[2020]{34K11}

\keywords{Second-order differential equation, Oscillation criteria, DEPCAG, functional differential equations, Hybrid dynamics.}

\date{}

\dedicatory{Dedicated to the memory of Prof. Manuel Pinto Jim\'enez.}
\maketitle
\begin{abstract}
This work presents two simple criteria for determining the oscillatory nature of solutions to second-order differential equations with deviated arguments. These criteria extend the (Leighton-Wintner)-type criteria established by G.Q. Wang and S.S. Cheng in \cite{Wang_cheng}, considering a generalized piecewise constant argument. 
Finally, we provide some examples.
\end{abstract}


\section{Introduction}
It is well known that oscillatory behavior frequently occurs in nature and can also involve piecewise-constant functions.
An example of this phenomenon is discussed in L. Dai’s book \cite{DAI}, where the author examines the oscillatory motion of a spring-mass system subject to piecewise constant forces. The system under study is given by  
\[
mx''(t) + cx'(t) + kx(t) = f(t, x([t])),
\]  
where \( f(t, x([t])) = Ax([t]) \) or \( f(t, x([t])) = B\cos(x([t])) \), and \( [\cdot] \) denotes the greatest integer function. A notable example of such a system is the Geneva wheel, a mechanism commonly used in watches (see also \cite{DAI_SINGH}).\\

\noindent In \cite{AK2}, as a generalization of A. Myshkis' work~\cite{203} on differential equations with deviating arguments, Marat Akhmet introduced an interesting class of differential equations of the form  
\begin{equation}
    y^{\prime}(t) = f(t, y(t), y(\gamma(t))), \label{depcag_eq}
\end{equation}
where \( \gamma(t) \) is a \emph{piecewise constant argument of generalized type}.  
The function \( \gamma(t) \) is defined as follows:  
Let \( \left( t_n \right)_{n\in\mathbb{Z}} \) and \( \left( \zeta_n \right)_{n\in\mathbb{Z}} \) be sequences satisfying  
\[
t_n < t_{n+1}, \quad \forall n \in \mathbb{Z},
\]  
with  
\[
\lim_{n\rightarrow\infty} t_n = \infty, \quad \lim_{n\rightarrow -\infty} t_n = -\infty,  
\]  
and let \( \zeta_n \in [t_n, t_{n+1}] \).  The function \( \gamma(t) \) is locally constant and then defined as  
\[
\gamma(t) = \zeta_n, \quad \text{for } t \in I_n = [t_n, t_{n+1}).
\]  
A fundamental example of such a function is  
$\gamma(t) = [t],$
where \( [\cdot] \) denotes the greatest integer function, which remains constant over each interval \( [n, n+1[ \) for \( n \in \mathbb{Z} \).

When a piecewise constant argument is introduced, the interval \( I_n \) can be decomposed into two subintervals: an advanced interval and a delayed interval, defined as  
\[
I_n = I_n^{+} \cup I_n^{-}, \quad \text{where} \quad I_n^{+} = [t_n, \zeta_n] \quad \text{and} \quad I_n^{-} = [\zeta_n, t_{n+1}].
\]  
Indeed, $$t\in I_k^+\implies t-\gamma(t)\leq 0,\qquad t\in I_k^-\implies t-\gamma(t)\geq 0.$$

Differential equations like \eqref{depcag_eq} are known as \emph{differential equations with piecewise constant argument of generalized type} (\emph{DEPCAG}).\\
One of their remarkable properties is that their solutions remain continuous despite the discontinuities of \( \gamma(t) \). Assuming the solutions of \eqref{depcag_eq} are continuous, integrating from \( t_n \) to \( t_{n+1} \) leads to an associated difference equation. As a result, these differential equations exhibit hybrid dynamics, incorporating both continuous and discrete characteristics (see \cite{AK2, P2011, Wi93}).  

For example, in \cite{Torres_proyecciones}, the author introduced the piecewise constant argument
\[
\gamma(t) = \left[\frac{t}{m}\right]m + \alpha m, \quad \text{where } m > 0 \text{ and } 0 \leq \alpha \leq 1.
\]
This definition implies that
\[
\left[ \frac{t}{m}\right]m + \alpha m = (n+\alpha )m, \quad \text{for } t \in I_n = [nm, (n+1)m).
\]
From these conditions, the advanced and delayed subintervals are determined by
\[
t- \gamma (t) \leq 0 \quad \Leftrightarrow \quad t \leq (n+\alpha)m, \quad \text{and} \quad t-\gamma (t) \geq 0 \quad \Leftrightarrow \quad t \geq (n+\alpha)m.
\]
Thus, the two subintervals can be written as
\[
I_n^{+} = [nm, (n+\alpha)m),\quad I_n^{-} = [(n+\alpha)m, (n+1)m).
\]

\begin{figure}[h!]
\centering
\includegraphics[scale=0.28]{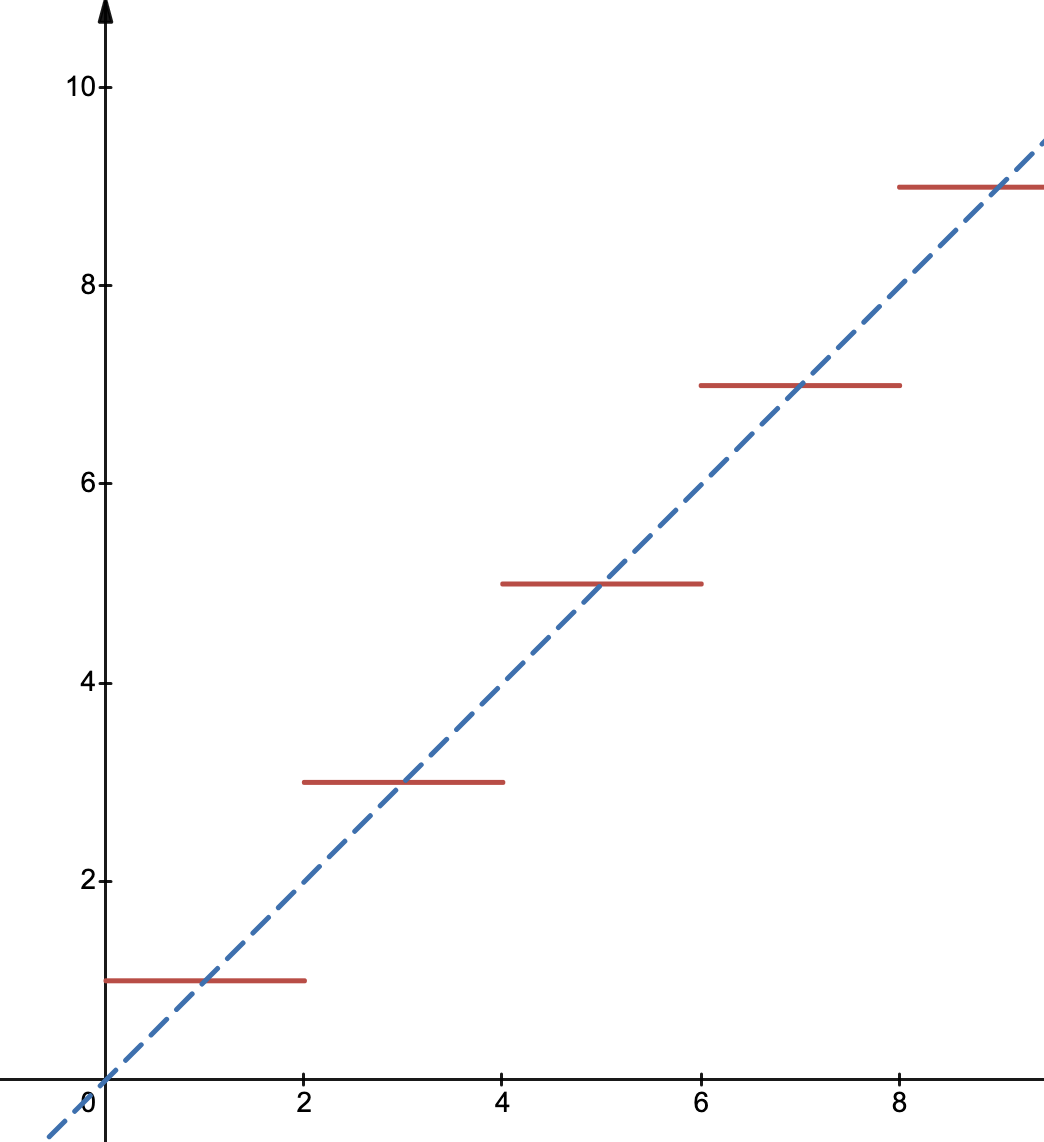}
\caption{$f(t)=2\left[\frac{t}{2}\right]+1$. An example of the previous piecewise constant argument, with $m=2$ and $\alpha=0.5$. }\label{fig1}
\end{figure}

\section{Recent Developments in Second-Order Differential Equations with Deviating Arguments}

In this section, we present some recent advancements in the study of second-order differential equations with deviating arguments. \\

\noindent In \cite{Yuan} (2003), R. Yuan investigated the existence of almost and quasi-periodic solutions for the following second-order differential equation with piecewise constant argument:
\[
    x''(t)+a(t)x(t)=\alpha x([t])+f(t),
\]
where \( a(t) \) is a 1-periodic continuous function, \( \alpha \neq 0 \), and \( f \) is a continuous function. The author also demonstrated that periodic and unbounded solutions can coexist in the equation
\[
    x''(t)+\omega^2x(t)=\alpha x([t])+f(t),
\]
which differs from the case \( \alpha=0 \). This phenomenon arises due to the piecewise constant argument and underscores a key distinction between ordinary differential equations and differential equations with piecewise constant arguments.\\

In \cite{Wang_Cheng_2006} (2006), G-Q. Wang and S.S. Cheng, utilizing Mawhin's continuation theorem, established the existence of periodic solutions for the second-order Rayleigh differential equation with piecewise constant argument:
\[
    x''(t)+f(t,x'(t))+g(t,x([t-k]))=0,
\]
where \( k\in\mathbb{Z}^+ \), and \( f(t,x) \) and \( g(t,x) \) are continuous on \( \mathbb{R}^2 \), satisfying \( f(t,0)=0 \) for all \( t\in\mathbb{R} \), and the periodicity conditions \( f(t+\omega,x)=f(t,x) \) and \( g(t+\omega,x)=g(t,x) \) for some \( \omega>0 \).\\

\noindent In \cite{Berek1} (2011), H. Bereketoglu, G. Seyhan, and F. Karakoc analyzed the second-order differential equation with piecewise constant mixed arguments:
\[
    x''(t)-a^2x(t)=bx([t-1])+cx([t])+dx([t+1]),
\]
where \( a, b, c, d \in \mathbb{R} \setminus \{0\} \). They proved the existence and uniqueness of solutions and established that the zero solution is a global attractor. Additionally, they explored the oscillatory behavior, non-oscillation properties, and periodicity of the solutions.\\

In \cite{Akhmet_mass_spring} (2023), M. Akhmet et al. examined a scalar undamped mass-spring system subject to piecewise constant forces of the form
\[
    mx''(t)+kx(t)=Ax(\gamma(t)),
\]
where \( A\in\mathbb{R} \), \( t\in\mathbb{R} \), and \( \gamma(t) \) is a generalized piecewise constant argument defined by \( \gamma(t)=t_k \) for \( t\in[t_k,t_{k+1}) \), \( k\in\mathbb{Z} \). The authors analyzed the solutions using the method of steps. \\

\noindent In \cite{fernandez_labora} (2024), S. Buedo-Fern\'andez, D. Cao Labora, and R. Rodr\'iguez-L\'opez studied the nonlinear second-order functional differential equation with piecewise constant arguments:
\[
    x''(t)=g(t,x(t),x'(t),x([t]),x'([t])),\quad t\in[0,T],
\]
subject to the boundary conditions \( x(0)=x(T) \) and \( x'(0)=x'(T)+\lambda \), where \( \lambda \in \mathbb{R} \), and \( g:[0,T]\times \mathbb{R}^4\to \mathbb{R} \) is continuous on \( [0,T] \setminus \{1,2,\ldots,[T]\} \times \mathbb{R}^4 \) and satisfies the conditions  
\[
    \lim_{t\to n^-}g(t,x,y,u,v),\quad \lim_{t\to n^+}g(t,x,y,u,v)=g(n,x,y,u,v)
\]
for finite limits. The authors established the existence of solutions within a certain region by approximation techniques.
This type of differential equation has applications in thermostat systems, where functional terms in the temperature and its rate of change at specific instants regulate the system's behavior.

\section{Aim of the work}
In \cite{Wang_cheng} (2004), \textit{Gen-Qiang Wang} and \textit{Sui Sun Cheng}
studied the following Second-order DEPCA 
\begin{equation}
\begin{tabular}{ll}
$(r(t)x(t)')^{\prime }+f(t,x([t]))=0$ & $t\geq 0,$ 
\end{tabular}
\label{depca_Wang_Cheng}
\end{equation}

\noindent Using certain integrability properties of the coefficients involved, the authors established an oscillatory Leighton-type criterion for \eqref{depca_Wang_Cheng}.\\
It is important to note that \( \gamma(t) = [t] \) is a particular case of a piecewise constant argument, where \( t_n = n = \zeta_n \) for \( n \in \mathbb{Z} \).\\

Inspired by \cite{leighton1950, leighton1952, Wang_cheng} and \cite{Wintner}, we establish two (Leighton-Wintner)-type oscillation criteria for the following DEPCAG:
\begin{equation}\label{osc_depcag_primera}
    (r(t)x(t)')^{\prime }+f(t,x(\gamma(t)))=0, \quad t\geq \tau,
\end{equation}
under certain hypotheses on the coefficients, as in \eqref{depca_Wang_Cheng}, but now with slight modifications due to the presence of a generalized piecewise constant argument \( \gamma(t) \).\\

Our work is structured as follows: first, we present some definitions and auxiliary results. Next, we state the main results. Finally, we provide examples that illustrate the effectiveness of our approach.

\section{Auxiliary Results}
\noindent During this work, we will use the following classical definitions of oscillation:
\begin{definition}~\par
\vskip1mm
A function $f(t)$ defined on $[t_0, \infty)$ is called oscillatory if there exist two sequences $(a_n),\,(b_n)\subset [t_0,\infty)$ such that $a_n\to\infty,\,b_n\to\infty$ as $n\to\infty$ and 
$f(a_n)\leq 0\leq f(b_n),\,\, \forall n\geq M,$ where $M$ is big enough. I.e., 
if the function does not eventually become strictly positive or negative, it is classified as oscillatory; otherwise, it is classified as non-oscillatory. 
\end{definition}
We remark that the last definition could also be interpreted in an asymptotic context.\\

Consider the following second-order DEPCAG:
\begin{equation}
\label{osc_DEPCAG}
(r(t)x'(t))' + f(t, x(\gamma(t))) = 0, \quad t \geq \tau,
\end{equation}
where $r(t)$ is a continuous and locally integrable function defined on $[\tau, \infty)$, and $f(t, x)$ is continuous on $[\tau, \infty) \times (-\infty, \infty)$, satisfying $x f(t, x) > 0$ for $t \geq \tau$ and $x \neq 0$. 

\noindent The function $\gamma(t)$ is a piecewise constant argument of generalized type such that $\gamma(t) = \zeta_n$ if $t \in I_n = [t_n, t_{n+1})$. 

\noindent Moreover, there exist locally integrable functions $p(t)$ and $\phi(x)$ such that $p(t)$ is continuous and nonnegative on $[\tau, \infty)$, and $\phi(x)$ is continuously differentiable and nondecreasing on $(-\infty, \infty)$, with $x \phi(x) > 0$ for $x \neq 0$, and
\[
f(t, x) \geq p(t) \phi(x), \quad x \neq 0,\ t \geq \tau.
\]

\begin{definition}{\cite{Wang_cheng,DAI,Yuan}}~\par
\vskip1mm
A continuous function $x(t)$ is a solution of \eqref{osc_DEPCAG} on $[\tau,\infty)$ if:
\begin{itemize}
\item [(i)] $x'(t)$ is continuously differentiable on $[\tau,\infty).$
\item[(i)] $(r(t)x'(t))'$ exists for all $t\in [\tau,\infty)$, except possibly at times $\{t_{k}\}_{k\in\mathbb{N}}$, where it has one-sided limits.
\item[(ii)] $x(t)$ satisfies \eqref{osc_DEPCAG} on the intervals of the form $[t_{k},t_{k+1}), k\in\mathbb{N}$.
\end{itemize}
\end{definition}

In the following, we will prove a useful Lemma that allows us to establish the main result:

\begin{lemma}\label{Wang_cheng_DEPCAG_lemma}
Let $x(t)$ be a solution of \eqref{osc_DEPCAG} such that there is some $M\geq 0$ and $x(t)\geq 0$ for $t\geq M$, $\gamma(t)$ be a generalized piecewise constant argument.\\
If 
\begin{equation}\label{condicion_lema1}
\displaystyle{\int_\tau^\infty \frac{1}{r(s)}ds=\infty},
\end{equation}
then $x'(t)\geq 0$ for all $\{t_k\}_{k\geq N},$ for some $N$ sufficiently large such that  $t_k\geq M$.
\end{lemma}
\begin{proof}
The lemma will be proved by contradiction. Suppose there exists \( M \in \mathbb{N} \) such that \( x'(t_k) < 0 \) for all \( t_k \geq M \). Let \( x'(t_k) = -\alpha \), with \( \alpha > 0 \).  

Since \( \gamma(t) = \zeta_k \) for all \( t \in I_k = [t_k, t_{k+1}) \), and given that \( p(t), r(t) > 0 \) and \( x\phi(x) > 0 \) for \( x \neq 0 \), it follows that \eqref{osc_DEPCAG} satisfies
\[
(r(t)x'(t))' = -f(t, x(\zeta_k)) \leq -p(t)\phi(x(\zeta_k)) \leq 0, \quad \forall \zeta_k \geq M.
\]
Hence, \( (r(t)x'(t))' \) is non-increasing on \( I_k \) for all \( t_k \geq M \).  

Therefore, for \( t \in I_k^- = [\zeta_k, t_{k+1}[ \), we have
\[
x'(t_{k+1}) \leq \frac{r(\zeta_k)}{r(t_{k+1})}x'(\zeta_k).
\]
Since \( (r(t)x'(t))' \) is non-increasing on \( I_k^- \), we also get
\[
x'(t) \leq \frac{r(\zeta_k)}{r(t)}x'(\zeta_k), \quad \forall t \in I_k^-.
\]
Integrating the above inequality yields
\[
x(t) \leq x(\zeta_k) + r(\zeta_k)x'(\zeta_k)\int_{\zeta_k}^{t} \frac{1}{r(s)}\,ds.
\]
Since \( x(t) \) is continuous, taking the limit as \( t \to t_{k+1} \) gives
\begin{equation} \label{parte1}
x(t_{k+1}) \leq x(\zeta_k) + r(\zeta_k)x'(\zeta_k)\int_{\zeta_k}^{t_{k+1}} \frac{1}{r(s)}\,ds.
\end{equation}
Proceeding similarly for \( t \in I_k^+ = [t_k, \zeta_k] \), we have
\begin{equation} \label{parte2}
x'(\zeta_k) \leq \frac{r(t_k)}{r(\zeta_k)}x'(t_k).
\end{equation}
Since \( (r(t)x'(t))' \) is non-increasing on \( I_k^+ \), we also have
\[
x'(t) \leq \frac{r(t_k)}{r(t)}x'(t_k), \quad \forall t \in I_k^+.
\]
Integrating the above inequality gives
\[
x(t) \leq x(t_k) + r(t_k)x'(t_k)\int_{t_k}^{t} \frac{1}{r(s)}\,ds.
\]
Taking the limit as \( t \to \zeta_k \), we obtain
\begin{equation} \label{parte3}
x(\zeta_k) \leq x(t_k) + r(t_k)x'(t_k)\int_{t_k}^{\zeta_k} \frac{1}{r(s)}\,ds.
\end{equation}

Now, applying \eqref{parte2} and \eqref{parte3} into \eqref{parte1}, we get
\[
x(t_{k+1}) \leq x(t_k) + r(t_k)x'(t_k)\int_{t_k}^{t_{k+1}} \frac{1}{r(s)}\,ds.
\]
Since \( x'(t_k) \leq -\alpha \), we have
\[
x(t_{k+1}) \leq x(t_k) - \alpha r(t_k)\int_{t_k}^{t_{k+1}} \frac{1}{r(s)}\,ds.
\]
Proceeding inductively, we obtain
\[
x(t_{k+n}) \leq x(t_k) - \alpha r(t_k)\int_{t_k}^{t_{k+n}} \frac{1}{r(s)}\,ds.
\]
According to condition \eqref{condicion_lema1}, the right-hand side tends to \( -\infty \) as \( n \to \infty \), which contradicts the assumption that \( x(t) \geq 0 \) for all \( t \geq M \).
\end{proof}

\begin{remark}
It is important to note that \( x(t_k) \) is defined in terms of \( x(\zeta_k) \), which is a critical detail. The advanced interval \( I_k^+ = [t_k, \zeta_k] \) and the delayed interval \( I_k^- = [\zeta_k, t_{k+1}[ \) must be considered in order to correctly define the solution of a DEPCAG over the interval \( [t_k, t_{k+1}] \) (see \cite{P2011}).
\end{remark}
\section{Main Results}
We are now in a position to prove our first oscillation criterion:

\begin{theorem}\label{wang_cheng_DEPCAG_criterio1}
    Assume that Lemma \ref{Wang_cheng_DEPCAG_lemma} holds, and suppose that
    \begin{equation}\label{wang_cheng_condicion_2}
        \int_\tau^\infty p(s) \, ds = \infty.
    \end{equation}
    Then, every solution of \eqref{osc_DEPCAG} is oscillatory. 
\end{theorem}

\begin{proof}
    Again, we will prove the theorem by contradiction. Suppose that \eqref{osc_DEPCAG} has a non-oscillatory solution \( x(t) \). We can also assume that \( x(t) \geq 0 \) for \( t \geq \tau \). By Lemma \ref{Wang_cheng_DEPCAG_lemma}, we know that \( x'(t) \geq 0 \) for \( t \geq \tau \).

    Define
    \begin{equation}
        w(t) = \dfrac{r(t)x'(t)}{\phi(x(\gamma(t)))}.
    \end{equation}
    We observe that \( w(t) \geq 0 \) for \( t \geq \tau \) and \( w(t_k^-) \geq 0 \) for all \( k \in \mathbb{N} \). By \eqref{osc_DEPCAG}, we have
    \begin{equation}\label{parte4}
        w'(t) = \dfrac{-f(t,x(\gamma(t)))}{\phi(x(\gamma(t)))} \leq -p(t), \quad t \in [t_k, t_{k+1}).
    \end{equation}

    Now, by the definition of the piecewise constant argument, 
    \[
        \lim_{t \to t_{k+1}^-} \gamma(t) = \zeta_k, \quad \gamma(t_{k+1}) = \zeta_{k+1},
    \]
    and since \( \phi(x) \) is non-decreasing and \( x\phi(x) > 0 \) for \( x \neq 0 \), we obtain
    \begin{equation}\label{desigualdad_tk1}
        w(t_{k+1}) = \dfrac{r(t_{k+1})x'(t_{k+1})}{\phi(x(\zeta_{k+1}))} \leq \dfrac{r(t_{k+1})x'(t_{k+1})}{\phi(x(\zeta_k))} = w(t_{k+1}^-).
    \end{equation}

    Next, we integrate \eqref{parte4} over \( I_k^- = [\zeta_k, t_{k+1}] \) and \( I_k^+ = [t_k, \zeta_k] \).

    \medskip
    \textbf{Step 1: Integrate over \( I_k^- \).} For \( t \in [\zeta_k, t_{k+1}) \), we get
    \begin{equation*}
        w(t) \leq w(\zeta_k) - \int_{\zeta_k}^{t} p(s)\,ds.
    \end{equation*}
    Taking the limit as \( t \to t_{k+1}^- \), we obtain
    \begin{equation}\label{parte5}
        w(t_{k+1}^-) \leq w(\zeta_k) - \int_{\zeta_k}^{t_{k+1}} p(s)\,ds.
    \end{equation}

    \medskip
    \textbf{Step 2: Integrate over \( I_k^+ \).} For \( t \in [t_k, \zeta_k) \), we have
    \begin{equation*}
        w(t) \leq w(t_k) - \int_{t_k}^{t} p(s)\,ds.
    \end{equation*}
    Taking the limit as \( t \to \zeta_k \), we get
    \begin{equation}\label{parte6}
        w(\zeta_k) \leq w(t_k) - \int_{t_k}^{\zeta_k} p(s)\,ds.
    \end{equation}

    Applying \eqref{parte6} in \eqref{parte5}, we obtain
    \begin{equation}\label{recursiva_criterio2}
        w(t_{k+1}^-) - w(t_k) \leq -\int_{t_k}^{t_{k+1}} p(s)\,ds.
    \end{equation}

    Now, using \eqref{parte4}, \eqref{desigualdad_tk1}, and \eqref{recursiva_criterio2}, we get
    \begin{align*}
        w(t_{k+2}^-) - w(t_k) &= w(t_{k+2}^-) - w(t_{k+1}) + w(t_{k+1}) - w(t_k) \\
        &\leq w(t_{k+2}^-) - w(t_{k+1}) + w(t_{k+1}^-) - w(t_k) \\
        &\leq -\int_{t_{k+1}}^{t_{k+2}} p(s)\,ds - \int_{t_k}^{t_{k+1}} p(s)\,ds \\
        &= -\int_{t_k}^{t_{k+2}} p(s)\,ds.
    \end{align*}

    Proceeding inductively, we obtain
    \begin{equation}\label{recursiva_criterio3}
        w(t_{k+n}^-) \leq w(t_k) - \int_{t_k}^{t_{k+n}} p(s)\,ds.
    \end{equation}

    Finally, as \( n \to \infty \), the right-hand side of the last expression tends to \( -\infty \), which contradicts the fact that \( w(t_k) \geq 0 \) for all \( t_k \geq M \). Hence, the solution \( x(t) \) is oscillatory.
\end{proof}

Next, we will present our second oscillatory criterion:
\begin{theorem}
    Suppose that \eqref{condicion_lema1} holds. Let $\varepsilon>0$. \\
    If 
    \begin{equation}\label{condicion_1_teo2}
        \displaystyle{\int_{\varepsilon}^\infty\frac{1}{\phi(u)}du<\infty, \quad \int_{-\varepsilon}^{-\infty}\frac{1}{\phi(u)}du<\infty}
    \end{equation}
    and
    \begin{equation}\label{condicion_2_teo2}
        \displaystyle{\int_{\tau}^\infty p(u)du<\infty, \quad 
        \int_{\tau
   }^{\infty}\frac{1}{r(s)}\left(\int_{t_{k(s)+1}}^\infty p(u)du\right)ds=\infty},
    \end{equation}
    where  $k(\tau)$ is the unique integer such that $\tau\in I_{k(\tau)}=[t_k,t_{k+1})$, then every solution of \eqref{osc_DEPCAG} is oscillatory.
\end{theorem}
\begin{proof}
    Suppose that \eqref{osc_DEPCAG} has no oscillatory solutions. Without loss of generality, assume that $x(t)\geq 0,\,\forall t\geq \tau.$ By Lemma \ref{Wang_cheng_DEPCAG_lemma}, we have $x'(t)\geq 0,$ $\forall t\geq \tau.$ Thus, $x(t)$ is non-decreasing on $[\tau,\infty).$ \\
    Consequently,  
    \begin{equation}\label{teo2_0}
    (r(t)x'(t))'=-f(t,x(\zeta_k))\leq -p(t)\phi(x(\zeta_k)),\quad \text{for all }t\in[t_k,t_{k+1}).
    \end{equation}
    Now, integrating \eqref{teo2_0} in $I_k^+=[t_k,\zeta_k],$ we get
    \begin{equation}\label{teo2_1}
    r(t)x'(t)-r(t_k)x'(t_k)\leq -\phi(x(\zeta_k))\int_{t_k}^t p(u)du.
    \end{equation}
By continuity of $x'(t)$, taking $t\to\zeta_k$ we see that
\begin{equation*}
    r(\zeta_k)x'(\zeta_k)-r(t_k)x'(t_k)\leq -\phi(x(\zeta_k))\int_{t_k}^{\zeta_k} p(u)du.
    \end{equation*}
Hence, we have
\begin{equation}\label{teo2_2}
    x'(t_k)\geq \frac{r(\zeta_k)}{r(t_k)}x'(\zeta_k)+\frac{\phi(x(\zeta_k))}{r(t_k)}\int_{t_k}^{\zeta_k} p(u)du.
\end{equation}
Next, integrating \eqref{teo2_0} in $I_k^-=[\zeta_k,t_{k+1}),$ we get
 \begin{equation*}
    r(t)x'(t)-r(\zeta_k)x'(\zeta_k)\leq -\phi(x(\zeta_k))\int_{\zeta_k}^t p(u)du.
    \end{equation*}
By the left continuity of $x'(t)$, taking $t\to t_{k+1}$ we see that
\begin{equation*}
    r(t_{k+1})x'(t_{k+1}^-)-r(\zeta_k)x'(\zeta_k)\leq -\phi(x(\zeta_k))\int_{\zeta_k}^{t_{k+1}} p(u)du.
    \end{equation*}
    From the last expression, we can conclude that
    \begin{equation}\label{teo2_3}
    x'(\zeta_k)\geq \frac{r(t_{k+1})}{r(\zeta_k)}x'(t^-_{k+1})+\frac{\phi(x(\zeta_k))}{r(\zeta_k)}\int_{\zeta_k}^{t_{k+1}} p(u)du.
\end{equation}
Applying \eqref{teo2_3} in \eqref{teo2_2}, we obtain
 \begin{equation}\label{teo2_4}
    x'(t_k)\geq \frac{r(t_{k+1})}{r(t_k)}x'(t^-_{k+1})+\frac{\phi(x(\zeta_k))}{r(t_k)}\int_{t_k}^{t_{k+1}} p(u)du.
\end{equation}
Now, considering $s,t$ such that $t_k\leq s\leq t\leq  \zeta_k=\gamma(s),$ integrating \eqref{teo2_0} and taking $t\to \zeta_k$, we get
\begin{equation*}
    r(\zeta_k)x'(\zeta_k)-r(s)x'(s)\leq -\phi(x(\zeta_k))\int_{s}^{\zeta_k} p(u)du.
    \end{equation*}
Then, from \eqref{teo2_3}, it follows that 
\begin{eqnarray*}
    x'(s)&\geq& \frac{r(\zeta_k)}{r(s)}x'(\zeta_k)+\frac{\phi(x(\zeta_k))}{r(s)}\int_{s}^{\zeta_k} p(u)du.\\
    &\geq&  \frac{r(\zeta_k)}{r(s)}x'(\zeta_k)\\
    &\geq& \frac{r(\zeta_k)}{r(s)}\left( \frac{r(t_{k+1})}{r(\zeta_k)}x'(t^-_{k+1})+\frac{\phi(x(\zeta_k))}{r(\zeta_k)}\int_{\zeta_k}^{t_{k+1}} p(u)du\right).
\end{eqnarray*}
That is,  
 \begin{equation}\label{teo2_5}
    x'(s)\geq \frac{r(t_{k+1})}{r(s)}x'(t^-_{k+1})+\frac{\phi(x(\zeta_k))}{r(s)}\int_{\gamma(s)}^{t_{k+1}} p(u)du.
\end{equation}
Also, from \eqref{teo2_4}, we obtain a lower bound for $x(t_{k+1}^-)$. Then, we have
\begin{eqnarray*}
    x'(s) &\geq& \frac{r(t_{k+1})}{r(s)}\left( \frac{r(t_{k+2})}{r(t_{k+1})}x'(t^-_{k+2})+\frac{\phi(x(\zeta_{k+1}))}{r(t_{k+1})}\int_{t_{k+1}}^{t_{k+2}} p(u)du\right)\\
    &&+\frac{\phi(x(\zeta_k))}{r(s)}\int_{\gamma(s)}^{t_{k+1}} p(u)du.
\end{eqnarray*}
Therefore, we see that
\begin{eqnarray}\label{teo2_6}
   x'(s)&\geq & \frac{r(t_{k+2})}{r(s)}x'(t^-_{k+2})+\frac{\phi(x(\zeta_{k+1}))}{r(s)}\int_{t_{k+1}}^{t_{k+2}} p(u)du\nonumber \\
   &&+
    \frac{\phi(x(\zeta_k))}{r(s)}\int_{\gamma(s)}^{t_{k+1}} p(u)du.
\end{eqnarray}
Again, by \eqref{teo2_4}, we have
\begin{equation*}
x'(t_{k+2})\geq  \frac{r(t_{k+3})}{r(t_{k+2})}x'(t^-_{k+3})+\frac{\phi(x(\zeta_{k+2}))}{r(t_{k+2})}\int_{t_{k+2}}^{t_{k+3}} p(u) du
\end{equation*}
Applying the last expression in \eqref{teo2_6} we obtain
\begin{eqnarray*}
   x'(s)&\geq & \frac{r(t_{k+2})}{r(s)}\left(\frac{r(t_{k+3})}{r(t_{k+2})}x'(t^-_{k+3})+\frac{\phi(x(\zeta_{k+2}))}{r(t_{k+2})}\int_{t_{k+2}}^{t_{k+3}} p(u)du\right)\\
   &&+\frac{\phi(x(\zeta_{k+1}))}{r(s)}\int_{t_{k+1}}^{t_{k+2}} p(u)du+
    \frac{\phi(x(\zeta_k))}{r(s)}\int_{\gamma(s)}^{t_{k+1}} p(u)du.
\end{eqnarray*}
I.e., 
\begin{eqnarray*}
   x'(s)&\geq & \frac{1}{r(s)}\left(r(t_{k+3})x'(t^-_{k+3})+\phi(x(\zeta_{k+2}))\int_{t_{k+2}}^{t_{k+3}} p(u)du \right.\\
   &&\left.+\phi(x(\zeta_{k+1}))\int_{t_{k+1}}^{t_{k+2}} p(u)du+\phi(x(\zeta_k))\int_{\gamma(s)}^{t_{k+1}} p(u)du\right).
\end{eqnarray*}
Hence, inductively, due to the positivity of the coefficients and $x'(t_j)\geq 0$, if $s\in I_{k(s)}^+$ and $t\in I_{n}=[t_{n},t_{n+1}),$ where $k(s)$ is the unique integer such that $s\in I_{k(s)}^+=[t_{k(s)},\zeta_{k(s)}]$, we have
\begin{eqnarray*}
   x'(s)&\geq & \frac{1}{r(s)}\left(\sum_{j=k(s)+1}^{k(s)+1+n}\phi(x(\zeta_{j}))\int_{t_{j}}^{t_{j+1}} p(u)du +\phi(x(\gamma(s)))\int_{\gamma(s)}^{t_{k+1}} p(u)du\right).
\end{eqnarray*}
Now, as $\phi(x)$ is non-decreasing, we have
\begin{eqnarray*}
   x'(s)\geq  \frac{\phi(x(s))}{r(s)}\int_{\gamma(s)}^{t_{k(s)+1+n}} p(u)du,
\end{eqnarray*}
or
\begin{eqnarray}\label{teo2_7}
   \frac{x'(s)}{\phi(x(s))}\geq  \frac{1}{r(s)}\int_{\gamma(s)}^{t_{k(s)+1+n}} p(u)du\geq \frac{1}{r(s)}\int_{t_{k(s)+1}}^{t_{k(s)+1+n}} p(u)du.
\end{eqnarray}
Taking $n\to\infty$, by \eqref{condicion_2_teo2}, we have
\begin{eqnarray}\label{teo2_8}
   \frac{x'(s)}{\phi(x(s))}\geq  \frac{1}{r(s)}\int_{t_{k(s)+1}}^{\infty} p(u)du,\quad \text{s }\in I_k^+=[t_k,\zeta_k].
\end{eqnarray}
Next, integrating \eqref{teo2_8} for $s\in [t_k,\zeta_k]$, we see that
\begin{eqnarray}
   \int_{t_k}^{\zeta_k}\frac{x'(s)}{\phi(x(s))}ds\geq  \int_{t_k}^{\zeta_k}\frac{1}{r(s)}\left(\int_{t_{k(s)+1}}^{\infty} p(u)du\right) ds. \label{cuentafinal1}
\end{eqnarray}
In the same way, integrating \eqref{teo2_8} for $s\in [\zeta_k,t_{k+1})$ we get
\begin{eqnarray}
   \int_{\zeta_k}^{t_{k+1}}\frac{x'(s)}{\phi(x(s))}ds\geq  \int_{\zeta_k}^{t_{k+1}}\frac{1}{r(s)}\left(\int_{t_{k(s)+1}}^{\infty} p(u)du\right) ds. \label{cuentafinal2}
\end{eqnarray}
Thus, from \eqref{cuentafinal1} and \eqref{cuentafinal2} we get
\begin{eqnarray*}
   \int_{t_k}^{t_{k+1}}\frac{x'(s)}{\phi(x(s))}ds\geq  \int_{t_k}^{t_{k+1}}\frac{1}{r(s)}\left(\int_{t_{k(s)+1}}^{\infty} p(u)du\right) ds.
\end{eqnarray*}
As $\phi(x),x'(s)>0,$ and using $z=x(s)$, we have
\begin{eqnarray*}
    \int_{x(t_k)}^{x(t_{k+1})}\frac{1}{\phi(z)}dz\geq \int_{t_k}^{t_{k+1}}\frac{1}{r(s)}\left(\int_{t_{k(s)+1}}^{\infty} p(u)du\right) ds.
\end{eqnarray*}
In this way, it is not difficult to see that
\begin{eqnarray*}
   \sum_{j=k(\tau)+1}^n\int_{x(t_j)}^{x(t_{j+1})}\frac{1}{\phi(z)}dz\geq \sum_{j=k(\tau)+1}^n\int_{t_j}^{t_{j+1}}\frac{1}{r(s)}\left(\int_{t_{k(s)+1}}^\infty p(u)du\right)ds. 
\end{eqnarray*}
Hence, if 
\begin{eqnarray*}
    \int_{x(\tau)}^{\infty}\frac{1}{\phi(z)}dz&&=\int_{x(\tau)}^{x(t_{k(\tau)+1})}\frac{1}{\phi(z)}dz+\sum_{j=k(\tau)+1}^\infty\int_{x(t_j)}^{x(t_{j+1})}\frac{1}{\phi(z)}dz,\\
\int_{\tau}^{\infty}\frac{1}{r(s)}\left(\int_{t_{k(s)+1}}^\infty p(u)du\right)ds&&=\int_{\tau}^{t_{k(\tau)+1}}\frac{1}{r(s)}\left(\int_{t_{k(s)+1}}^\infty p(u)du\right)ds\\
   &&+\sum_{j=k(\tau)+1}^\infty\int_{t_j}^{t_{j+1}}\frac{1}{r(s)}\left(\int_{t_{k(s)+1}}^\infty p(u)du\right)ds,
\end{eqnarray*}
by \eqref{condicion_1_teo2}, we obtain
\begin{equation*}
   \int_{x(\tau)}^\infty \frac{1}{\phi(z)}dz\geq 
   \int_{\tau
   }^{\infty}\frac{1}{r(s)}\left(\int_{t_{k(s)+1}}^\infty p(u)du\right)ds,
\end{equation*} 
where $x(\tau)$ is some initial condition for \eqref{osc_DEPCAG}.\\
Finally, this yields a contradiction, since the left-hand side is finite while, by \eqref{condicion_2_teo2}, the right-hand side is infinite.
\end{proof}

\begin{remark}
\begin{itemize}
\item[] 
\item Notably, our results hold regardless of the choice of \( \gamma(t) \).
\item If we consider $\gamma(t) = [t]$, then $t_{k(s)+1} = [s] + 1$, thereby recovering the result obtained by Wang and Cheng in \cite{Wang_cheng}, where the authors used the condition  
\[
\int_{\tau}^{\infty} \frac{1}{r(s)} \left( \int_{[s]+1}^{\infty} p(u) \, du \right) ds = \infty.
\]

\end{itemize}
\end{remark}

\section{Examples}
In this section, we will give some examples that show the applicability of our results.

\begin{example}
Consider the following DEPCAG:
\[
x''(t) + kx(\gamma(t)) = 0,
\]
where $k > 0$ and $\gamma(t)$ is a generalized piecewise constant argument.

\noindent Since $f(t, x) = k\phi(x) = kx$ and $r(t) = 1$, by Theorem~\ref{wang_cheng_DEPCAG_criterio1}, all solutions are oscillatory.

In fact, if we consider $\gamma(t) = [t]$, $k = 2$, $x'(0) = 0$, and $x(0) = 1$, then the discrete solution of
\begin{equation}\label{ejemplo1}
x''(t) + 2x([t]) = 0
\end{equation}
is\footnote{Calculated at \url{https://www.wolframcloud.com}} 
\[
x(n) = \dfrac{\left(\frac{1}{2}(1 - i\sqrt{7})\right)^n(\sqrt{7} - i) + \left(\frac{1}{2}(1 + i\sqrt{7})\right)^n(\sqrt{7} + i)}{2\sqrt{7}},\quad n \in \mathbb{N} \cup \{0\},
\]
which is oscillatory.

\begin{figure}[h!]
\centering
\includegraphics[scale=0.23]{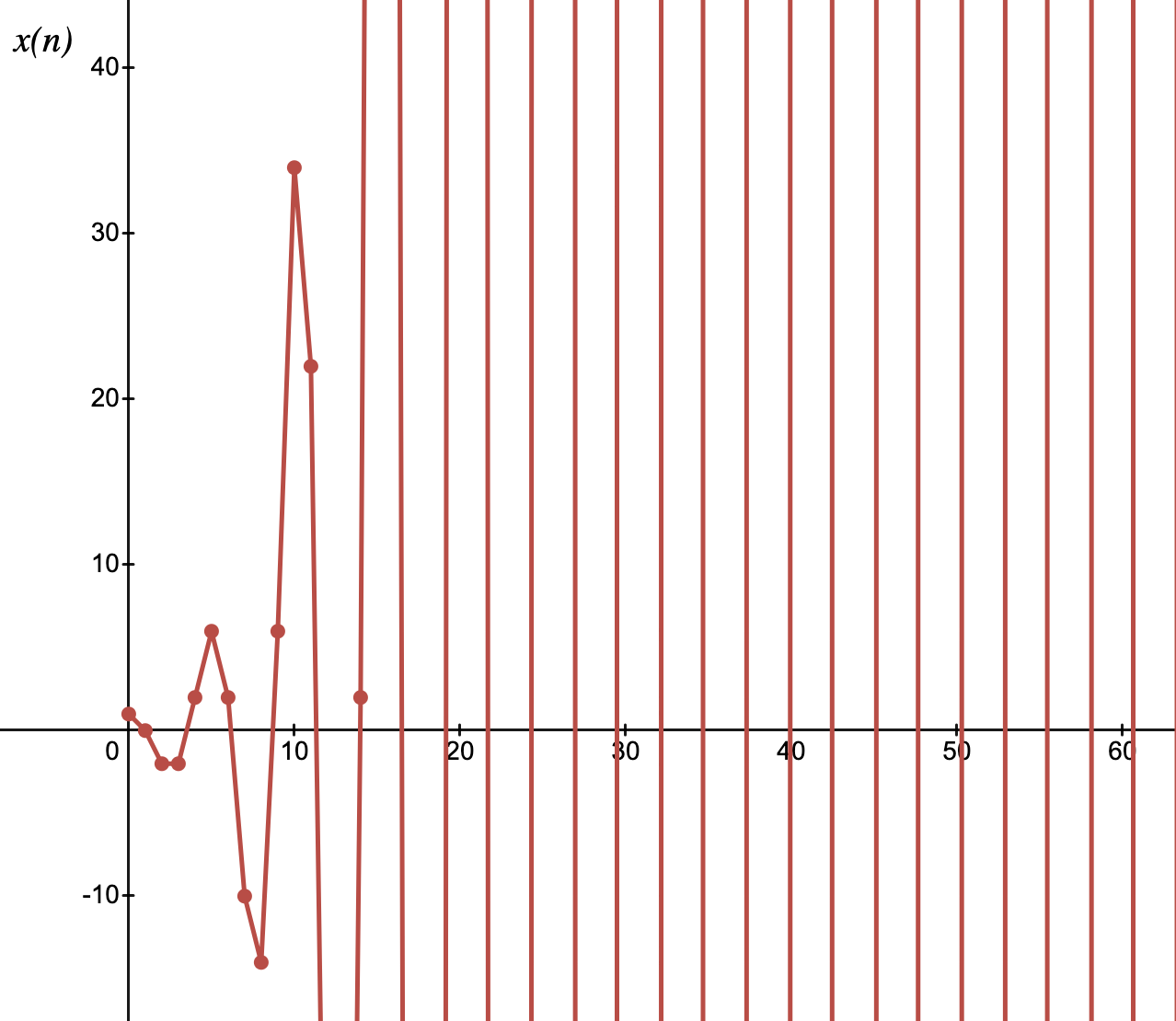}
\caption{Discrete solution $x(n)$ of \eqref{ejemplo1}, with $x(0)=1$ and $x'(0)=0$.}\label{fig2}
\end{figure}
\end{example}
\newpage
\begin{example}
Taking into account the example given in \cite{Wang_cheng}, consider the following DEPCAG:
\begin{equation}\label{general}
(\exp(-t)x'(t))' + x(\gamma(t))\exp(t^2 + (x(\gamma(t)))^2) = 0, \quad t \geq 0,
\end{equation}
where $\gamma(t)$ is any piecewise constant argument, $r(t) = \exp(-t)$, $\phi(x) = x$, $p(t) = \exp(t^2)$, and $f(t,x) = x\exp(t^2 + x^2)$.

\noindent It is not difficult to see that
\[
f(t,x) = x\exp(t^2 + x^2) \geq x\exp(t^2) = p(t)\phi(x), \quad x \neq 0,\ t \geq 0,
\]
and that $r(t)$ and $p(t)$ are such that conditions \eqref{condicion_lema1} and \eqref{wang_cheng_condicion_2} are satisfied. Then, by Theorem~\ref{wang_cheng_DEPCAG_criterio1}, all solutions of \eqref{general} are oscillatory.

\end{example}
\section*{Acknowledgements}
\noindent This work is lovingly dedicated to the memory of Prof.\ Manuel Pinto Jim\'enez, whose profound influence, mentorship, and unwavering support continue to inspire. Lastly, sincere thanks are extended to DESMOS PBC for granting permission to use the images created with the DESMOS graphing calculator, accessible at  
\url{https://www.desmos.com/calculator}.

\end{document}